\documentclass[a4paper,12pt]{amsart}
\usepackage{amsfonts}
\usepackage{amssymb}
\usepackage{ifthen}
\usepackage{graphics}
\usepackage{amsmath}
\usepackage{latexcad}
\usepackage{graphics,graphicx}
\include{graphix}

\nonstopmode \numberwithin{equation}{section}
\newtheorem{thm}{Theorem}
\newtheorem{cor}{Corollary}
\newtheorem{lem}{Lemma}
\newtheorem{prop}{Proposition}
\newtheorem{opbl}{Open Problem}

\newtheorem{claim}{Claim}
\newtheorem{conj}[equation]{Conjecture}

\theoremstyle{definition}
\newtheorem{defn}{Definition}
\newtheorem{case}{Case}
\newtheorem{examp}[equation]{Example}
\newtheorem{prob}[equation]{Problem}
\newtheorem{ques}[equation]{Question}
\newtheorem{rem}{Remark}

\newcounter {own}
\def\theown {\thesection       .\arabic{own}}

\newenvironment{pf}[1][]{%
 \vskip 3mm
 \noindent
 \ifthenelse{\equal{#1}{}}%
  {{\slshape Proof. }}%
  {{\slshape #1.} }%
 }%
{\qed\bigskip}

\newcounter{alphabet}
\newcounter{tmp}
\newenvironment{Thm}[1][]{\refstepcounter{alphabet}%
\bigskip%
\noindent%
{\bf Theorem \Alph{alphabet}}%
\ifthenelse{\equal{#1}{}}{}{ (#1)}%
{\bf .} \itshape}{\vskip 8pt}

\makeatletter   
\newcommand{\Ref}[1]{\@ifundefined{r@#1}{}{\setcounter{tmp}{\ref{#1}}\Alph{tmp}}}
\makeatother
\newenvironment{Lem}[1][]{\refstepcounter{alphabet}%
\bigskip%
\noindent%
{\bf Lemma \Alph{alphabet}}%
{\bf .} \itshape}{\vskip 8pt}

\newcommand{\IR}{{\mathbb R}}

\def\be{\begin{equation}}
\def\ee{\end{equation}}

\newcommand{\bee}{\begin{enumerate}}
\newcommand{\eee}{\end{enumerate}}

\newcommand{\blem}{\begin{lem}}
\newcommand{\elem}{\end{lem}}
\newcommand{\bthm}{\begin{thm}}
\newcommand{\ethm}{\end{thm}}
\newcommand{\bcor}{\begin{cor}}
\newcommand{\ecor}{\end{cor}}
\newcommand{\beg}{\begin{examp}}
\newcommand{\eeg}{\end{examp}}
\newcommand{\begs}{\begin{examples}}
\newcommand{\eegs}{\end{examples}}
\newcommand{\bdefe}{\begin{defn}}
\newcommand{\edefe}{\end{defn}}
\newcommand{\bprob}{\begin{prob}}
\newcommand{\eprob}{\end{prob}}
\newcommand{\bques}{\begin{ques}}
\newcommand{\eques}{\end{ques}}
\newcommand{\bei}{\begin{itemize}}
\newcommand{\eei}{\end{itemize}}
\newcommand{\bop}{\begin{opbl}}
\newcommand{\eop}{\end{opbl}}
\newcommand{\bcl}{\begin{claim}}
\newcommand{\ecl}{\end{claim}}

\newcommand{\bca}{\begin{case}}
\newcommand{\eca}{\end{case}}
\newcommand{\bcon}{\begin{conj}}
\newcommand{\econ}{\end{conj}}
\newcommand{\bcons}{\begin{conjs}}
\newcommand{\econs}{\end{conjs}}
\newcommand{\bprop}{\begin{prop}}
\newcommand{\eprop}{\end{prop}}
\newcommand{\br}{\begin{rem}}
\newcommand{\er}{\end{rem}}
\newcommand{\brs}{\begin{rems}}
\newcommand{\ers}{\end{rems}}
\newcommand{\bo}{\begin{obser}}
\newcommand{\eo}{\end{obser}}
\newcommand{\bos}{\begin{obsers}}
\newcommand{\eos}{\end{obsers}}
\newcommand{\bpf}{\begin{pf}}
\newcommand{\epf}{\end{pf}}
\newcommand{\ba}{\begin{array}}
\newcommand{\ea}{\end{array}}
\newcommand{\beq}{\begin{eqnarray}}
\newcommand{\beqq}{\begin{eqnarray*}}
\newcommand{\eeq}{\end{eqnarray}}
\newcommand{\eeqq}{\end{eqnarray*}}

\newcommand{\ds}{\displaystyle}

\newcounter{minutes}
\divide\time by 60
\newcounter{hours}
\multiply\time by 60 \addtocounter{minutes}{-\time}

\begin{document}

\bibliographystyle{amsplain}
\title{Freely quasiconformal maps and distance ratio metric}
\thanks{$^\dagger$ File:~\jobname .tex,
          printed: \number\year-\number\month-\number\day,
          \thehours.\ifnum\theminutes<10{0}\fi\theminutes}
\author{Yaxiang  Li}
\address{Y. Li,  College of Science,
Central South University of
Forestry and Technology, Changsha,  Hunan 410004, People's Republic
of China} \email{yaxiangli@163.com}

\author{Saminathan Ponnusamy  $^\dagger $}
\address{S. Ponnusamy,
Indian Statistical Institute (ISI), Chennai Centre,
SETS (Society for Electronic Transactions and security),
MGR Knowledge City, CIT Campus, Taramani,
Chennai 600 113, India.
}
\email{samy@isichennai.res.in, samy@iitm.ac.in}

\author{Matti Vuorinen}
\address{Matti Vuorinen Department of Mathematics, University of Turku, FIN-20014 Turku, Finland.}
\email{vuorinen@utu.fi}

\date{}
\keywords{ Quasihyperbolic metric, quasiconformal mapping, FQC mapping, distance ratio metric.\\
$
^\dagger$ {\tt This authors is on leave from  the Department of Mathematics,
Indian Institute of Technology Madras, Chennai-600 036, India}}

\begin{abstract}
Suppose that $E$ and $E'$ denote real Banach spaces with dimension at least $2$ and that $D\subset E$ and $D'\subset E'$ are
domains. In this paper, we establish, in terms  of the $j_D$ metric,
a necessary and sufficient condition for the homeomorphism $f: E \to E'$ to be FQC. Moreover, we give, in terms  of the $j_D$ metric, a sufficient condition for the homeomorphism $f: D\to D'$  to be FQC. On the other hand, we show that this condition is not necessary.

\end{abstract}

\maketitle\pagestyle{myheadings} \markboth{Y. Li, S. Ponnusamy, and M. Vuorinen}{Freely quasiconformal maps and distance ratio metric}

\section{Introduction and main results}\label{sec-1}

During the past few decades, modern mapping theory and the geometric
theory of quasiconformal maps has been studied from several points of view. These studies include Heinonen's work on metric measure spaces \cite{Hei}, Koskela's study
of maps with finite distortion \cite{KKM} and V\"ais\"al\"a's work about quasiconformality in
infinite dimensional Banach spaces \cite{Vai6-0, Vai6, Vai6', Vai4, Vai8}.
 Our study is motivated by  V\"ais\"al\"a's theory of freely quasiconformal maps in the setup of Banach spaces \cite{Vai6-0, Vai6, Vai6'}.  The basic tools in V\"ais\"al\"a's theory are metrics and the notion of
 uniform continuity between metric spaces, in particular the norm metric, the quasihyperbolic metric and the distance ratio metric are used.
 We begin with some basic definitions and the statements of our results.

Throughout the paper, we always assume that $E$ and $E'$ denote real Banach
spaces with dimension at least $2$, and that $D\subset E$ and $D'\subset E'$ are
domains. The norm of a vector $z$ in
$E$ is written as $|z|$, and for each pair of points $z_1$, $z_2$ in $E$,
the distance between them is denoted by $|z_1-z_2|$. The
distance from $z\in D$ to the boundary $\partial D$ of $D$ is denoted by $d_D(z)$.

For each pair of points $z_1$, $z_2$ in $D$, the {\it distance ratio
metric} $j_D(z_1,z_2)$ between $z_1$ and $z_2$ is defined by
$$j_D(z_1,z_2)=\log\left(1+\frac{|z_1-z_2|}{\min\{d_D(z_1),d_D(z_2)\}}\right).$$

The {\it quasihyperbolic length} of a rectifiable arc or a path
$\alpha$ in the norm metric in $D$ is the number (cf.
\cite{GP,Vai6-0}):
$$\ell_k(\alpha)=\int_{\alpha}\frac{|dz|}{d_{D}(z)}.
$$

For each pair of points $z_1$, $z_2$ in $D$, the {\it quasihyperbolic distance}
$k_D(z_1,z_2)$ between $z_1$ and $z_2$ is defined in the usual way:
$$k_D(z_1,z_2)=\inf\ell_k(\alpha),
$$
where the infimum is taken over all rectifiable arcs $\alpha$
joining $z_1$ to $z_2$ in $D$. Gehring and Palka \cite{GP} introduced the quasihyperbolic metric of
a domain in $\IR^n$. Many of the basic properties of this metric may
be found in \cite{Geo}. We remark that the quasihyperbolic metric has
been recently studied by many people (cf. \cite{HIMPS, HPWW, Kle,klen,krt, RT,rt2}).


\bdefe \label{def1.7} Let $D\not=E$ and $D'\not=E'$ be metric
spaces, and let $\varphi:[0,\infty)\to [0,\infty)$ be a growth
function, that is, a homeomorphism with $\varphi(t)\geq t$. We say
that a homeomorphism $f: D\to D'$ is {\it $\varphi$-semisolid} if
$$k_{D'}(f(x),f(y))\leq \varphi(k_D(x,y))$$
for all $x$, $y\in D$, and {\it $\varphi$-solid} if both $f$ and $f^{-1}$
satisfy this condition.\edefe

The special case $\varphi(t)=Mt \;(M>1)$ gives the $M$-quasihyperbolic maps or briefly $M$-QH. More precisely, $f$ is called $M$-QH if
$$\frac{k_D(x,y)}{M}\leq k_{D'}(f(x),f(y))\leq Mk_D(x,y)$$ for  all  $x$ and $y$ in $D$.

We say that $f$ is {\it fully $\varphi$-semisolid}
(resp. {\it fully $\varphi$-solid}) if $f$ is
$\varphi$-semisolid (resp. $\varphi$-solid) on every  subdomain of $D$. In particular,
when $D=E$, the subdomains are taken to be proper ones in $D$. Fully $\varphi$-solid mappings are also called {\it freely
$\varphi$-quasiconformal mappings}, or briefly $\varphi$-FQC mappings.

If $E=E'=\IR^n$, then $f$ is FQC if and only if $f$ is
quasiconformal (cf. \cite{Vai6-0}). See \cite{Avv,Vai1, Mvo1} for definitions and
properties of $K$-quasiconformal mappings, or briefly $K$-QC mappings.

It is well-known that for all $z_1$, $z_2$ in $D$, we have
(cf. \cite{Vai6-0})
\begin{equation}\label{eq(0000)}j_D(z_1, z_2)\leq \inf \log\left(1+\frac{\ell(\alpha)}{\min\{d_D(z_1),d_D(z_2)\}}\right)\leq\inf\ell_k(\alpha)=k_{D}(z_1, z_2),\end{equation} where $\alpha$ is a rectifiable arc joining $z_1$ and $z_2$ in $D$.
Hence, in the study of FQC maps, it is natural to ask whether we could use $j_D$ metric to describe FQC or not. In fact, we get the following conditions for a homeomorphism to be FQC.

\begin{thm}\label{thm2.0} A homeomorphism $f: E \to  E'$ is $\varphi_1$-{\rm FQC} if and only if for every proper subdomain $D\subset E$, we have
\beq\label{eq(0)}\varphi_2^{-1}(j_{D}(x,y))\leq j_{D'}(f(x),f(y))\leq \varphi_2(j_{D}(x,y))\eeq
with $x,y\in D,$ where $\varphi_1$ and $\varphi_2$ are self-homeomorphisms of $[0,\infty)$ with $\varphi_i(t)\geq t$ $(i=1,2)$ for all $t$, and $\varphi_1$, $\varphi_2$ depend only on each other. \end{thm}

\begin{thm}\label{thm3.0}Let $\varphi:[0,\infty)\to [0,\infty)$ be a homeomorphism with $\varphi(t)\geq t$ for all $t$, and let $f: D\varsubsetneq E \to D'\varsubsetneq E'$ be a homeomorphism. If for every subdomain $D_1\subset D$, we have
\beq\label{eq(00)}\varphi^{-1}(j_{D_1}(x,y))\leq j_{D'_1}(f(x),f(y))\leq \varphi(j_{D_1}(x,y))\eeq
for all $x,y\in D_1,$ then $f$ is $\varphi_1$-{\rm FQC} with $\varphi_1=\varphi_1(\varphi).$  Moreover, if $D_1=D$ and $\varphi(t)=Mt$ ($M\geq 1$), then $f$ is $M_1$-{\rm QH} with $M_1=M_1(M)$.\end{thm}

\begin{thm}\label{exm1}The converse of Theorem \ref{thm3.0} is not true.\end{thm}

%

%

 The proofs of Theorems \ref{thm2.0}, \ref{thm3.0} and \ref{exm1} will be given in Section \ref{sec-4}.
 In Section \ref{sec-2}, some necessary preliminaries will be introduced.

\section{Preliminaries}\label{sec-2}

 For an open ball with center $x\in E$ and radius $r> 0$ we use the notation $\mathbb{B}(x,r)$. The boundary of the ball is
denoted by $\mathbb{S}(x,r)$.  The closed
line segment with endpoints $z_1\in E$ and $z_2\in E$ is denoted by $[z_1, z_2]$.

For convenience, in what follows, we always assume that $x$, $y$, $z$, $\ldots$
denote points in $D$ and $x'$, $y'$, $z'$, $\ldots$
the images in $D'$ of $x$, $y$, $z$, $\ldots$
under $f$, respectively. Also we assume that $\alpha$, $\beta$, $\gamma$, $\ldots$
denote curves in $D$ and $\alpha'$, $\beta'$, $\gamma'$, $\ldots$  the images in $D'$ of
$\alpha$, $\beta$, $\gamma$, $\ldots$
under $f$, respectively.

\bdefe \label{def1.3} A domain $D$ in $E$ is called $c$-{\it
uniform} in the norm metric provided there exists a constant $c$
with the property that each pair of points $z_{1},z_{2}$ in $D$ can
be joined by a rectifiable arc $\alpha$ in $ D$ satisfying (cf. \cite{MS,Vai, Vai6})

\bee
\item\label{wx-4} $\ds\min_{j=1,2}\ell (\alpha [z_j, z])\leq c\,d_{D}(z)
$ for all $z\in \alpha$, and

\item\label{wx-5} $\ell(\alpha)\leq c\,|z_{1}-z_{2}|$,
\eee

\noindent where $\ell(\alpha)$ denotes the length of $\alpha$,
$\alpha[z_{j},z]$ the part of $\alpha$ between $z_{j}$ and $z$.
\edefe

It is well known that each ball $B\subset E$ is $2$-uniform and a half space $H\subset E$ is $c$-uniform for all $c>2$ (cf. \cite[Example 10.4]{Vai8}).

In \cite{Vai6},  V\"ais\"al\"a characterized uniform domains by the quasihyperbolic metric.

\begin{Thm}\label{thm0.1} {\rm (\cite[Theorem 6.16]{Vai6})}
For a domain $D$ in $E$, the following are quantitatively equivalent: \bee

\item $D$ is a $c$-uniform domain;
\item $k_D(z_1,z_2)\leq c'\,
 \log\left(1+\frac{\ds{|z_1-z_2|}}{\ds\min\{d_{D}(z_1),d_{D}(z_2)\}}\right)$ for every pair of points $z_1$, $z_2\in D$;
\item $k_D(z_1,z_2)\leq c'_1\,
 \log\left(1+\frac{\ds{|z_1-z_2|}}{\ds\min\{d_{D}(z_1),d_{D}(z_2)\}}\right)+d$ for every pair of points $z_1$, $z_2\in D$.\eee
\end{Thm}

 In the case of domains in $ {\mathbb R}^n,$  the equivalence
  of items (1) and (3) in Theorem D is due to Gehring and Osgood \cite{Geo} and the
  equivalence of items (2) and (3) due to Vuorinen \cite{Mvo2}. 

\begin{Thm}\label{ThmO} $($\cite[Lemma 6.7]{Vai6}$)$~~Suppose that $D$ is a $c$-uniform domain and that $x_0\in D$. Then $D_0=D\setminus\{x_0\}$ is $c_0$-uniform
with $c_0=c_0(c)$.\end{Thm}

\begin{Thm}\label{ThmD} $($\cite[Lemma 6.5]{Vai6}$)$~~For $x_0\in E$, $E\setminus\{x_0\}$ is $10$-uniform.\end{Thm}
For the case of $G=\mathbb{R}^n\setminus\{0\}$ the sharp constant $d$ in the inequality $k_G(x,y)\le d j_G(x,y)$ was found by H. Linden \cite{li}, and it is $\pi/\log 3$.

%
%
%
%

 Recall that an arc $\alpha$ from $z_1$ to
$z_2$ is a {\it quasihyperbolic geodesic} if
$\ell_k(\alpha)=k_D(z_1,z_2)$. Each subarc of a quasihyperbolic
geodesic is obviously a quasihyperbolic geodesic. It is known that a
quasihyperbolic geodesic between every pair of points in $E$ exists if the
dimension of $E$ is finite, see \cite[Lemma 1]{Geo}. This is not
true in arbitrary spaces (cf. \cite[Example 2.9]{Vai4}).
In order to remedy this shortcoming, V\"ais\"al\"a introduced the following concepts \cite{Vai6}.

\bdefe \label{def1.4}Let $D\neq E$ and $c\geq 1$. An arc $\alpha\subset D$ is a $c$-neargeodesic if and only if $\ell_k(\alpha[x,y])\leq c\,k_D(x,y)$
for all $x, y\in \alpha$.
\edefe

In \cite{Vai6}, V\"ais\"al\"a proved the following property concerning
the existence of neargeodesics in $E$.

\begin{Thm}\label{LemA} $($\cite[Theorem 3.3]{Vai6}$)$
Let $\{z_1, z_2\}\subset D$ and $c>1$. Then there is a
$c$-neargeodesic in $D$ joining $z_1$ and $z_2$.
\end{Thm}

%
A map $f:X\to Y$ is {\it uniformly continuous} \cite[Section 3.1]{Vai6-0} if and only if there is $t_0\in (0,\infty]$ and an embedding $\varphi:[0,t_0)\to [0,\infty)$
with $\varphi(0)=0$ such that $$|f(x)-f(y)|\leq \varphi(|x-y|)$$ whenever $x,y\in X$ and $|x-y|\leq t_0$. We then say that $f$ is {\it $(\varphi, t_0)$-uniformly continuous}. If $t_0=\infty$, we briefly say that $f$ is {\it $\varphi$-uniformly continuous}.

\begin{Lem}\label{Thmp} $($\cite[Lemma 3.2]{Vai6-0}$)$~~Let $X$ be $c$-quasiconvex and let $f:X\to Y$ be a map. Then the following conditions are quantitatively equivalent:\bee
\item $f$ is $(\varphi, t_0)$-uniformly continuous;

\item $f$ is $\varphi$-uniformly continuous.

\eee\end{Lem}

%
%

In \cite{Vai6-0},  V\"ais\"al\"a proved the following theorem.

\begin{Thm}\label{ThmQ} $($\cite[Theorem 5.13]{Vai6-0}$)$~~Suppose that $f: E\to D'\subset E'$ is fully $\varphi$-semisolid. Then  \bee
\item $D'=E'$;

\item $f$ is $\phi$-{\rm FQC} with $\psi=\psi_{\varphi}$;

\item $f$ is $\eta$-{\rm QS} with $\eta=\eta_{\varphi}$.
\eee\end{Thm}

\section{The proofs of Theorems \ref{thm2.0}, \ref{thm3.0} and  \ref{exm1} }\label{sec-4}

Before the proofs of our main results, we list a series of lemmas.
From the proof of Theorem $5.7$ and summary $5.11$ in \cite{Vai6-0}, we get the following.

\begin{Lem}\label{LemO}Let $\varphi:[0,\infty)\to [0,\infty)$ be a homeomorphism with $\varphi(t)\geq t$ for all $t$. Suppose that $f:D\to D'$ is a homeomorphism. If for every point $x\in D$, $f:D\setminus\{x\}\to D'\setminus\{x'\}$ is $\varphi$-solid, then $f$ is $\psi$-{\rm FQC} with $\psi=\psi(\varphi).$\end{Lem}

\begin{Lem}\label{LemR} $($\cite[Lemma 6.25]{Vai6-0}$)$~~Suppose that $D\not=E$, $D'\not=E'$ and that $f: D\to D'$ is a $\theta$-{\rm QM} homeomorphism $($for definition see \cite {Vai2, Vai6-0}$)$. Then $$j_{D'}(x,y)\leq M j_D(x,y)+C$$  for all $x,y\in D$ with $M$ and $C$ depending only on $\theta$. \end{Lem}

The following result is from \cite{Vai2}.

\begin{Lem}\label{LemQ} $($\cite[Theorem 3.2]{Vai2}$)$~~Suppose that $f:X\to Y$ is $\eta$-{\rm QS}. Then $f$ is $\theta$-{\rm QM} where $\theta$ depends only on $\eta$. \end{Lem}

\begin{lem}\label{lem-2} Suppose that $f:E\to E'$ is $\varphi$-{\rm FQC}. Then for every proper subdomain $D\subset E$,
$$j_{D'}(x,y)\leq c j_D(x,y)+d$$ with $x,y \in D$,  where $c$, $d$ are constants depending only on $\varphi$.
\end{lem}
\bpf
By Theorem \Ref{ThmQ}, we know that $f$ is $\eta$-{\rm QS} with $\eta$ depending only on $\varphi$. Then for every proper subdomain $D\subset E$, $f:D\to D'$ is an $\eta$-{\rm QS} homeomorphism. Hence, Lemma \ref{lem-2} follows from Lemmas \Ref{LemR} and \Ref{LemQ}.
\epf

\begin{lem}\label{lem-1} Let $D\subset E$ be a domain, and let $x, y\in D$. If $|x-y|\leq sd_D(x)$ for some $s\in (0,1)$, then
$$k_{B_x}(x,y)\leq\frac{1}{1-s}\log\Big(1+\frac{|x-y|}{d_D(x)}\Big),$$ where $B_x=\mathbb{B}(x,d_D(x))$.
\end{lem}

\bpf For each $w\in [x,y]$, we have $d_{B_x}(w)\geq d_D(x)-|x-w|$. Then
\begin{eqnarray*}k_{B_x}(x,y)&\leq& \int_{[x,y]}\frac{|dw|}{d_{B_x}(w)}\leq \int^{d_D(x)}_{d_D(x)-|x-y|}\frac{dt}{t}
\\ \nonumber&=& \log\left(1+\frac{|x-y|}{d_D(x)-|x-y|}\right)\leq
\log\left(1+\frac{|x-y|}{(1-s)d_D(x)}\right)\\ \nonumber &\leq&\frac{1}{1-s}\log\left(1+\frac{|x-y|}{d_D(x)}\right),\end{eqnarray*}
from which our lemma follows.\epf

We remark that when $E=\IR^n$, Lemma \ref{lem-1} coincides with Lemma $3.7$  in \cite{Mvo1}.\medskip

Now we are ready to prove our main results.
\subsection{The proof of Theorem \ref{thm2.0}}
We first prove the sufficient part and we assume that $f:E\to E'$ is $\varphi_1$-FQC. Let $D$ be a proper subdomain in $E$. For $x,y\in D$, by symmetry, to prove \eqref{eq(0)}, we only need to prove the right hand side inequality. Without loss of generality, we may assume that $d_D(x)\leq d_D(y)$.

If $|x-y|\leq \frac{1}{2}d_D(x)$, let $B_x=\mathbb{B}(x,d_D(x))$. Then \eqref{eq(0000)} and Lemma \ref{lem-1} show that
\beq\label{theorem1-eq1} j_{D'}(x',y')&\leq&  k_{D'}(x',y')\leq \varphi_1(k_{D}(x,y))\\ \nonumber &\leq& \varphi_1(k_{B_x}(x,y))\leq \varphi_1\left(2 \log\left(1+\frac{|x-y|}{d_D(x)}\right)\right)\\ \nonumber &\leq& \varphi_1(2 j_{D}(x,y)).
\eeq

If $|x-y|> \frac{1}{2}d_D(x)$, then $$j_D(x,y)=\log\left(1+\frac{|x-y|}{d_D(x)}\right)> \log (3/2),$$ which, together with Lemma \ref{lem-2}, shows that
\beq\label{theorem1-eq2}j_{D'}(x',y')\leq c j_D(x,y)+d\leq \left(c+d/\log (3/2)\right)j_D(x,y),\eeq where $c$ and $d$ are constants depending only on $\varphi_1.$

Hence, by \eqref{theorem1-eq1} and \eqref{theorem1-eq2}, inequality \eqref{eq(0)} holds with $\varphi_2(t)=\max \{\varphi_1(2t), \big(c+d/\log (3/2)\big)t\}$.

In the following, we prove the necessary part. For fixed $a\in E$, let $D_a=E\setminus \{a\}$. Then from Theorem \Ref{ThmD} that both $D_a$ and $D_a'$ are $c$-uniform with $c=10$. Hence, Theorem \Ref{thm0.1} yields that there is a universal constant $c'$ such that for each pair of points $x,y\in D_a$,
$$k_{D_a'}(x',y')\leq c'j_{D_a'}(x',y')\;\; {\rm and} \;\;k_{D_a}(x,y)\leq c'j_{D_a}(x,y).$$
Therefore, it follows from \eqref{eq(0000)} and \eqref{eq(0)} that $$\varphi^{-1}_2\left(\frac{1}{c'}k_{D_a}(x,y)\right)\leq k_{D_a'}(x',y')\leq c'\varphi_2(k_{D_a}(x,y)),$$ which shows that $f$ is $c'\varphi_2$-solid in $D_a$, and so from Lemma \Ref{LemO} that $f:E\to E'$ is $\varphi_1$-FQC with $\varphi_1$ depending only on $\varphi_2$.

Hence the proof of Theorem \ref{thm2.0} is complete.
\qed

\subsection{The proof of Theorem \ref{thm3.0}}
We first prove the first part. Let $D_1$ be a subdomain of $D$. 
  By symmetry, we only need to prove that there exists a homeomorphism $\varphi_0$ such that for each pair of points $x,y \in D_1$, \beq\label{maineq}k_{D_1'}(x',y')\leq \varphi_0(k_{D_1}(x,y)).\eeq
Choose $0<t_0<1$ such that $\varphi(t_0)\leq \log(3/2)$. Let $x,y\in D_1$ be points with $k_{D_1}(x,y)\leq t_0$. Then \eqref{eq(00)} gives $$j_{D'_1}(x',y')\leq \varphi(j_{D_1}(x,y))\leq \varphi(k_{D_1}(x,y))\leq \log(3/2),$$ which shows that  $$|x'-y'|\leq \frac{1}{2}d_{D_1'}(x').$$ Let $B_{x'}=\mathbb{B}(x', d_{D_1'}(x'))$.
 Then by  Lemma \ref{lem-1}, \eqref{eq(0000)}and \eqref{eq(00)} , we obtain
\begin{eqnarray*}k_{D_1'}(x',y')&\leq& k_{B_{x'}}(x',y')\leq 2 \log\left(1+\frac{|x'-y'|}{d_{D_1'}(x')}\right)
\\ \nonumber &\leq&2 j_{D_1'}(x',y')\leq2 \varphi(j_{D_1}(x,y))\\ \nonumber &\leq&2\varphi(k_{D_1}(x,y)).\end{eqnarray*}
Hence, we see from Lemma  \Ref{Thmp} that \eqref{maineq} holds.

Now, we prove the second part. Let $a=1-e^{-\frac{1}{3M}}$. We only need to prove that there exists a constant $M_1=M_1(M)$ such that for each pair of points $x$, $y\in D$,   \beq\label{maineq1}k_{D'}(x',y')\leq M_1 k_{D}(x,y).\eeq We divide the proof into two cases.

 \bca\label{case1} $|x-y|\leq a d_D(x).$\eca Then $$d_D(y)\geq d_D(x)-|x-y|\geq (1-a)d_D(x).$$ Hence,
$$j_{D'}(x',y')\leq M j_{D}(x,y)= M \log\left(1+\frac{|x-y|}{\min\{d_{D}(x),d_D(y)\}}\right)\leq \frac{1}{3},
$$ which shows that $$|x'-y'|< \frac{1}{2}d_{D'}(x').$$ Let $B_x'=\mathbb{B}(x', d_{D_1'}(x'))$. Then Lemma \ref{lem-1} yields that
\begin{eqnarray*}k_{D'}(x',y')&\leq& k_{B_{x'}}(x',y')\leq 2 \log\left(1+\frac{|x'-y'|}{d_{D'}(x')}\right)
\\ \nonumber &\leq&2 j_{D'}(x',y')\leq2M j_{D_1}(x,y)\\ \nonumber &\leq&2Mk_{D}(x,y).\end{eqnarray*}
\bca\label{case2} $|x-y|\geq a d_D(x)$.\eca It follows from Theorem \Ref{LemA} that there exists a
$2$-neargeodesic $\gamma $ in $D$ joining $x$ and $y$. Let $x=z_{1}$, and let $z_{2}$ be the first intersection point of $\gamma$
with $\mathbb{S}(z_1,a d_{D}(z_1))$ in the direction from
$x$ to $y$. We let $z_{3}$ be the first intersection point of
$\gamma$ with $\mathbb{S}(z_2,a d_{D}(z_2))$ in the
direction from $z_{2}$ to $y$. By repeating this procedure, we get a
set $\{z_i\}_{i=1}^p$ of points
 in $\gamma$ such that $y$ is contained
in $\overline{\mathbb{B}}(z_{p},a d_{D}(z_p))$, but not in
$\overline{\mathbb{B}}(z_{p-1},a d_{D}(z_{p-1}))$. Obviously,
$p> 1$.
Hence Case \ref{case1} yields
\begin{eqnarray*}
k_{D'}(x', y')&\leq& \sum_{i=1}^{p-1} k_{D'}(z'_i, z'_{i+1})+k_{D}(z'_p, y')
\\ \nonumber
&\leq& 2M\left(\sum_{i=1}^{p-1} \;k_{D}(z_i,z_{i+1})+k_{D}(z_p,
y)\right)
\\ \nonumber
&\leq& 2M\ell_k(\gamma[x,y])
\\ \nonumber
&\leq&4M\,k_{D}(x,y).
\end{eqnarray*}
Hence, inequality \eqref{maineq1} holds, which complete the proof of Theorem \ref{thm3.0}.\qed

\subsection{The proof of  Theorem \ref{exm1}}

We prove this theorem by presenting two examples.
\begin{examp}\label{exm1'}Let $E=\IR^2\cong \mathbb{C}$ and $f$ be a conformal mapping of the unit disk $\mathbb{B}(0,1)=\{z: |z|<1\} (=D)$ onto $D'=\mathbb{B}(0,1)\setminus[0,1)$. There exist points $x,y\in D$ such that \eqref{eq(00)} does not hold.
 \end{examp}

{\bf Proof.} By \cite[Theorem 3]{Geo}, we know that conformal mapping is $M$-QH mapping for some constant $M\geq1$. Hence, $f$ is $\varphi$-FQC with $\varphi(t)=M t$.

Now, we prove that there is some domain $D_1\subset D$ such that \eqref{eq(00)} fails to hold.  Let $w_0'=(-\frac{1}{2},0)\in D'$ and let $D_1=D\setminus \{w_0\}$. Since each ball is $2$-uniform,  Theorem \Ref{ThmO} yields that $D_1$ is $c_1$-uniform, where $c_1$ is an universal constant. For each pair of points $x,y\in D_1$, we have
\beq\label{exam-eq}\frac{k_D(x,y)}{M}\leq k_{D'}(x',y')\leq Mk_D(x,y).\eeq

\begin{figure}
\begin{center}

\unitlength=1.1mm
\begin{picture}(76,42)
\thinlines
\drawthickdot{16.0}{22.0}
\thicklines
\drawdotline{16.0}{22.0}{28.0}{22.0}
\drawpath{60.0}{22.0}{72.0}{22.0}
\thinlines
\drawthickdot{60.0}{22.0}
\thicklines
\path(22.0,28.0)(22.0,28.0)(22.29,28.22)(22.59,28.45)(22.9,28.68)(23.2,28.9)(23.5,29.11)(23.81,29.33)(24.11,29.52)(24.43,29.72)
\path(24.43,29.72)(24.74,29.91)(25.04,30.09)(25.36,30.27)(25.66,30.45)(25.99,30.61)(26.29,30.77)(26.61,30.93)(26.93,31.09)(27.25,31.24)
\thinlines
\drawthickdot{20.0}{24.0}
\drawthickdot{14.0}{30.0}
\drawthickdot{12.0}{16.0}
\thicklines
\path(27.25,31.24)(27.56,31.38)(27.88,31.5)(28.2,31.63)(28.52,31.77)(28.84,31.88)(29.15,32.0)(29.49,32.11)(29.81,32.2)(30.13,32.31)
\thinlines
\drawthickdot{66.0}{24.0}
\drawthickdot{66.0}{20.0}
\thicklines
\drawcircle{60.0}{22.0}{24.0}{}
\drawcircle{16.0}{22.0}{24.0}{}
\path(56.0,22.0)(56.0,22.0)(56.0,22.06)(56.0,22.15)(56.0,22.22)(56.01,22.31)(56.01,22.38)(56.02,22.45)(56.04,22.52)(56.05,22.59)
\path(56.05,22.59)(56.08,22.66)(56.08,22.74)(56.12,22.79)(56.13,22.86)(56.16,22.93)(56.19,23.0)(56.22,23.06)(56.25,23.11)(56.27,23.18)
\path(56.27,23.18)(56.31,23.24)(56.36,23.29)(56.4,23.36)(56.44,23.4)(56.48,23.45)(56.51,23.52)(56.56,23.56)(56.62,23.61)(56.66,23.66)
\path(56.66,23.66)(56.72,23.72)(56.77,23.75)(56.83,23.81)(56.9,23.86)(56.95,23.9)(57.01,23.93)(57.08,23.97)(57.15,24.02)(57.22,24.06)
\path(57.22,24.06)(57.29,24.09)(57.36,24.13)(57.44,24.16)(57.51,24.2)(57.58,24.22)(57.68,24.27)(57.76,24.29)(57.83,24.33)(57.93,24.34)
\path(57.93,24.34)(58.01,24.38)(58.11,24.4)(58.19,24.43)(58.3,24.45)(58.4,24.47)(58.5,24.5)(58.59,24.5)(58.69,24.52)(58.8,24.54)
\path(58.8,24.54)(58.91,24.56)(59.01,24.58)(59.12,24.59)(59.23,24.61)(59.36,24.61)(59.48,24.63)(59.58,24.63)(59.72,24.63)(59.83,24.65)
\path(59.83,24.65)(59.95,24.65)(60.08,24.65)(60.22,24.65)(60.34,24.65)(60.48,24.65)(60.62,24.65)(60.76,24.65)(60.88,24.65)(61.04,24.65)
\path(61.04,24.65)(61.18,24.63)(61.31,24.63)(61.47,24.63)(61.62,24.61)(61.76,24.61)(61.91,24.59)(62.08,24.58)(62.23,24.56)(62.38,24.56)
\path(65.2,20.08)(65.59,20.04)(65.98,20.0)(66.0,20.0)
\path(65.29,20.08)(65.63,20.04)(65.98,20.0)(66.0,20.0)
\path(62.38,24.56)(62.55,24.54)(62.72,24.52)(62.87,24.5)(63.05,24.47)(63.22,24.45)(63.38,24.43)(63.55,24.4)(63.73,24.38)(63.91,24.36)
\path(63.91,24.36)(64.08,24.34)(64.27,24.31)(64.45,24.27)(64.63,24.25)(64.83,24.2)(65.01,24.18)(65.2,24.15)(65.4,24.11)(65.59,24.06)
\path(65.59,24.06)(65.8,24.02)(65.98,24.0)(66.0,24.0)
\path(56.0,22.0)(56.0,22.0)(56.0,21.91)(56.0,21.84)(56.0,21.75)(56.01,21.68)(56.01,21.61)(56.02,21.54)(56.04,21.45)(56.05,21.38)
\path(56.05,21.38)(56.08,21.31)(56.08,21.25)(56.12,21.18)(56.13,21.11)(56.16,21.06)(56.19,20.99)(56.22,20.93)(56.25,20.86)(56.27,20.81)
\path(56.27,20.81)(56.31,20.75)(56.36,20.68)(56.4,20.63)(56.44,20.58)(56.48,20.52)(56.51,20.47)(56.56,20.41)(56.62,20.36)(56.66,20.31)
\path(56.66,20.31)(56.72,20.27)(56.77,20.22)(56.83,20.18)(56.9,20.13)(56.95,20.09)(57.01,20.04)(57.08,20.0)(57.15,19.97)(57.22,19.93)
\path(57.22,19.93)(57.29,19.88)(57.36,19.86)(57.44,19.81)(57.51,19.79)(57.58,19.75)(57.68,19.72)(57.76,19.68)(57.83,19.65)(57.93,19.63)
\path(57.93,19.63)(58.01,19.61)(58.11,19.58)(58.19,19.56)(58.3,19.54)(58.4,19.52)(58.5,19.5)(58.59,19.47)(58.69,19.45)(58.8,19.43)
\path(58.8,19.43)(58.91,19.41)(59.01,19.4)(59.12,19.4)(59.23,19.38)(59.36,19.36)(59.48,19.36)(59.58,19.36)(59.72,19.34)(59.83,19.34)
\path(59.83,19.34)(59.95,19.34)(60.08,19.33)(60.22,19.33)(60.34,19.33)(60.48,19.33)(60.62,19.33)(60.76,19.33)(60.88,19.34)(61.04,19.34)
\path(61.04,19.34)(61.18,19.34)(61.31,19.34)(61.47,19.36)(61.62,19.36)(61.76,19.38)(61.91,19.38)(62.08,19.4)(62.23,19.41)(62.38,19.43)
\path(62.38,19.43)(62.55,19.45)(62.72,19.47)(62.87,19.49)(63.05,19.5)(63.22,19.52)(63.38,19.54)(63.55,19.58)(63.73,19.59)(63.91,19.63)
\path(63.91,19.63)(64.08,19.65)(64.27,19.68)(64.45,19.7)(64.63,19.74)(64.83,19.77)(65.01,19.81)(65.2,19.84)(65.4,19.88)(65.59,19.91)
\path(65.59,19.91)(65.8,19.95)(65.98,20.0)(66.0,20.0)
\thinlines
\drawthickdot{52.0}{22.0}
\thicklines
\path(30.13,32.31)(30.45,32.4)(30.77,32.49)(31.09,32.58)(31.43,32.65)(31.75,32.72)(32.08,32.79)(32.4,32.86)(32.72,32.91)(33.06,32.97)
\path(33.06,32.97)(33.38,33.02)(33.7,33.06)(34.04,33.11)(34.36,33.15)(34.7,33.18)(35.02,33.2)(35.36,33.22)(35.68,33.25)(36.02,33.25)
\path(36.02,33.25)(36.34,33.27)(36.68,33.27)(37.0,33.27)(37.34,33.27)(37.66,33.25)(37.99,33.25)(38.31,33.22)(38.65,33.2)(38.97,33.18)
\path(38.97,33.18)(39.31,33.15)(39.63,33.11)(39.97,33.08)(40.3,33.04)(40.62,32.99)(40.95,32.93)(41.29,32.88)(41.62,32.83)(41.94,32.77)
\path(41.94,32.77)(42.27,32.7)(42.59,32.63)(42.93,32.56)(43.26,32.49)(43.58,32.4)(43.91,32.33)(44.23,32.24)(44.55,32.15)(44.88,32.06)
\path(44.88,32.06)(45.2,31.95)(45.54,31.86)(45.86,31.75)(46.18,31.65)(46.5,31.54)(46.83,31.43)(47.15,31.31)(47.47,31.18)(47.79,31.06)
\path(47.79,31.06)(48.11,30.93)(48.41,30.81)(48.73,30.68)(49.05,30.54)(49.37,30.4)(49.69,30.25)(50.0,30.11)(50.31,29.97)(50.62,29.81)
\drawcenteredtext{38.0}{36.0}{\bf $f$}
\drawcenteredtext{23.0}{24.0}{$w_0$}
\drawvector{50.0}{30.0}{6.0}{2}{-1}
\thinlines
\drawcenteredtext{14.0}{28.0}{$x_0$}
\drawcenteredtext{9.5}{16.0}{$y_0$}
\drawcenteredtext{29.5}{22.0}{1}
\drawcenteredtext{52.0}{19.0}{$w_0'$}
\drawcenteredtext{62.0}{27.0}{$\gamma '$}
\drawcenteredtext{68.50}{25.0}{$x_0'$}
\drawcenteredtext{68.5}{19.0}{$y_0'$}
\drawcenteredtext{16.0}{6.0}{$D=\mathbb{B}(0, 1)$}
\drawcenteredtext{60.0}{6.0}{$D'=\mathbb{B}(0, 1) \setminus [0, 1)$}
\end{picture}

\end{center}
\caption{$\gamma' $ is a quasihyperbolic geodesic joining $x_0'$ and $y_0'$ in $D'$ \label{examfig1}}
\end{figure}
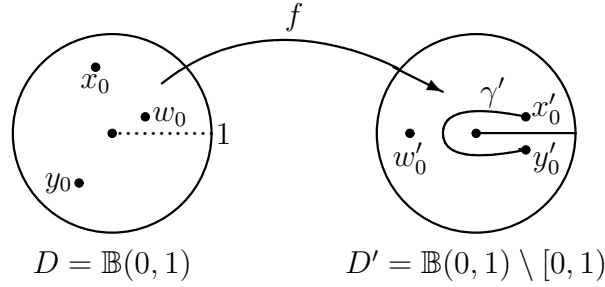

Take $x_0'$, $y_0'\in D_1'$ with $x_0'=(\frac{1}{2},t)$ and $y_0'=(\frac{1}{2},-t)$ (see Figure \ref{examfig1}).  Let $\gamma' $ be a quasihyperbolic geodesic joining $x'$ and $y'$ in $D'$. Then by \eqref{eq(0000)},
$$k_{D'}(x_0',y_0')\geq \log\left(1+\frac{\ell(\gamma')}{d_{D'}(x_0')}\right)\geq \log(1+\frac{1}{t}),$$
and
$$j_{D_1'}(x_0',y_0')=\log\left(1+\frac{|x_0'-y_0'|}{\min\{d_{D_1'}(x_0'),d_{D_1'}(y_0')\}}\right)=\log 3.$$
But \eqref{exam-eq} yields that
\begin{eqnarray*}j_{D_1}(x_0,y_0)&\geq& \frac{1}{c_1}k_{D_1}(x_0,y_0)\geq \frac{1}{c_1}k_{D}(x_0,y_0)\geq \frac{1}{c_1M}k_{D'}(x_0',y_0')\\ \nonumber &\geq&\frac{1}{c_1M}\log\left(1+\frac{1}{t}\right)
\rightarrow\infty, \;\; {\rm as}\;\; t\rightarrow0.
\end{eqnarray*}

Hence, there is no self-homeomorphism $\varphi$ of $[0,\infty)$ such that \eqref{eq(00)} hold. The proof of  Example \ref{exm1'} is complete.

Note that we could also choose $D_1=D$, and conclude from a similar proof that \eqref{eq(00)} does not hold.
\qed

\begin{examp}\label{exm2} Let $E$ be an infinite-dimensional separable real Hilbert space with an orthonormal basis $(e_j)_{j\in \mathbb{Z}}$ where all vectors $e_j$ have a constant norm. Setting $\gamma'_j=[e_{j-1},e_j]$ we obtain the infinite broken line $\gamma'=\cup\{\gamma'_j:j\in \mathbb{Z}\}$. Let $\gamma$ denote the line spanned by $e_1$, and let $D=\gamma+\mathbb{B}(r)$ with $r\leq \frac{1}{10}$ and $f$ be a locally $M$-bilipschitz homeomorphism from $D$ onto a neighborhood $D'$ of $\gamma'$ (for a detailed explanation we refer to see \cite[section 4.12]{Vai6-0} or \cite{vai2004}). Then there exist points $x,y\in D$ such that \eqref{eq(00)} does not hold. \end{examp}

{\bf Proof. }By \cite[Theorem 4.8]{Vai6-0}, we obtain that $f$ is $M^2$-QH. Let $x,y\in D$ with $x=\sqrt{2}e_1$, $y=m\sqrt{2}e_1$. Then $d_D(x)=d_D(y)=r$. Because $f$ is locally $M$-bilipschitz, we get $$d_{D'}(x')\geq \frac{r}{M}\;\; {\rm and}\;\; d_{D'}(y')\geq \frac{r}{M}.$$
Since $D'\subset \mathbb{B}(0,2)$, we find that
 $$j_{D'}(x',y')=\log\left(1+\frac{|x'-y'|}{\min\{d_{D'}(x'),d_{D'}(y')\}}\right)\leq \log\left(1+\frac{4M}{r}\right),$$
 but $$j_D(x,y)=\log\left(1+\frac{|x-y|}{\min\{d_{D}(x),d_{D}(y)\}}\right)=\log\left(1+\frac{\sqrt{2}(m-1)}{r}\right)\rightarrow\infty , $$ as $m\rightarrow\infty$.
 Hence \eqref{eq(00)} does not hold.\qed

\bigskip
\bigskip
{\bf Acknowledgements.} The first author was partially supported by a grant from Simons Foundation and Talent Introduction Foundation of Central South University of Forestry and Technology (No. 2013RJ005).
This research was completed during the visit of the first author to Indian Statistical Institute (ISI) Chennai Centre,
India. She thanks ISI Chennai Centre for the hospitality and the other supports.

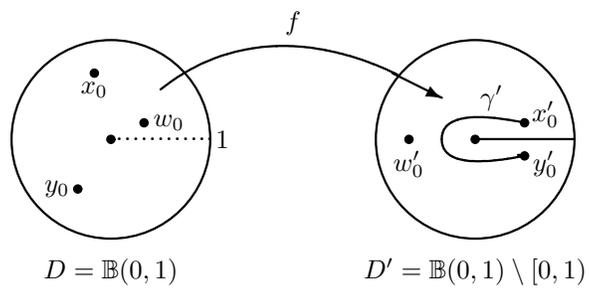
\end{document}